\documentclass[reqno,12pt]{amsart}
\usepackage{amsfonts}
\usepackage{amssymb}
\usepackage{bbm}
\usepackage{times}
\usepackage{bbm}
\usepackage{amsmath}
\usepackage{txfonts}
\usepackage{stmaryrd}
\usepackage{mathrsfs}

\font\bbold=msbm10

\def\DD{\hbox{\bbold D}}

\title{The Beurling-type theorem in the Bergman space $A^2_\alpha(\DD)$ for any $-1<\alpha<+\infty$
}

\author{Junfeng Liu}

\address{Junfeng Liu\\ Department of Mathematics, College of Science, Zhejiang University of Science and Technology, Hangzhou, 310023, P. R. China}   \email{jfliu997@sina.com}

\date{}
\newtheorem{theorem}{Theorem}

\begin{document}

\maketitle

\begin{abstract}
In this paper, we show that the Beurling-type theorem is true for a large class invariant subspaces of the shift operator $M_{z^N}$
of  multiplicity $N$ on the Bergman space $A^2_\alpha(\DD)$.
More specifically, every reducing subspace $H$ for the shift operator $S=M_{z^N}$ of multiplicity $N$ on the Bergman space $A^2_\alpha(\DD)$
has the property $H=[H\ominus SH]_{S,A^2_\alpha\left(\DD\right)}$ for any positive integer $N$ and any $-1<\alpha<+\infty$.
\end{abstract}



\section{Introduction}

It is well known that the structural theory of invariant subspaces for the shift operator on
Bergman spaces is so complicated that very little is known about it.
In this area, the most famous result is due to A. Aleman, S. Richter and C. Sundberg \cite{asc96}, who shown that
\textbf{\underline{every}} invariant subspace $H$ for the shift operator $S=M_z=M_{z^1}$ on the Bergman space $A^2(\DD)$ has the property
$$H=[H\ominus SH]_{S,A^2\left(\DD\right)},$$
where $[H\ominus SH]_{S,A^2\left(\DD\right)}$ denotes the smallest invariant subspace for the shift
operator $S$ on $A^2(\DD)$ which contains the set $H\ominus SH$.
In this case, it is often said that the Beurling-type theorem is true for the shift operator $S$ on the
Bergman spaces $A^2(\DD)~(=A^2_0(\DD))$.
Later S. Shimorn in \cite{sms01} and \cite{sms02} shown that the Beurling-type theorem
is true for the shift operator $S=M_z$ on the the Bergman space $A^2_\alpha(\DD)$ whenever $-1<\alpha <0$ and $A^2_\alpha(\DD)$
whenever $0<\alpha\leq 1$ respectively.

On the other hand, S. Shimorin \cite{sms01} conjectured that the Beurling-type theorem may fail for the shift operator $S=M_z$ on
$A^2_\alpha(\DD)$ with $\alpha>4$ (\cite{sms01} p.186); S. Shimorin \cite{sms02} stated that for $\alpha>1$, it seems that such a theorem
fails (\cite{sms02} p.1779).
Ulteriorly, by \cite{hp04}, the Beurling-type theorem is \textbf{\underline{not true}} for the shift operator $S=M_z$ on the
Bergman space $A^2_\alpha(\DD)$ when ever \textbf{$\underline{\alpha\geq 1.04}$}.
On the other words, it is not true that \textbf{\underline{every}} invariant subspace $H$ for the shift operator
$S=M_z$ on  $A^2_\alpha(\DD)$ has the property $H=[H\ominus SH]_{S,A^2\left(\DD\right)}$
whenever $\alpha\geq 1.04$.

In this paper, we show that the Beurling-type theorem is true for a large class invariant subspaces of the shift operator $M_{z^N}$
of  multiplicity $N$ on the Bergman space $A^2_\alpha(\DD)$.
More specifically, every reducing subspace $H$ for the shift operator $S=M_{z^N}$ of multiplicity $N (N=1,2,\cdot\cdot\cdot)$ on
the Bergman space $A^2_\alpha(\DD)$ $(-1<\alpha < +\infty)$
has the property $H=[H\ominus zH]_{S,A^2_\alpha\left(\DD\right)}$.

\section{preliminary}

For convenience, we first recall some basic concepts and facts from the references, which are used in this paper.

Let $\DD$ denote the open unit disk of the complex plane.
Let $-1<\alpha <+\infty$ and let $A^2_\alpha(\DD)$ be the (standard weighted) Bergman space defined by
$$A^2_\alpha(\DD)=\{f: f(z)=\sum^\infty_{n=0}a_nz^n,\ \ z\in \DD,\ \ \sum^\infty_{n=0}\omega_n |a_n|^2<+\infty\},$$
with the norm
$$\|f\|_{A^2_\alpha\left(\DD\right)}=\left(\sum^\infty_{n=0}\omega_n |a_n|^2\right)^{\frac{1}{2}},\ \ f(z)=\sum^\infty_{n=0}a_nz^n \in A^2_\alpha(\DD),$$
and the inner product

$$\langle f,g\rangle=\sum^\infty_{n=0}\omega_n a_n\overline{b_n}, \ \ f(z)=\sum^\infty_{n=0}a_nz^n \in A^2_\alpha(\DD),\ \
g(z)=\sum^\infty_{n=0}b_nz^n \in A^2_\alpha(\DD),$$
where

$$\omega_n=\frac{n!\Gamma(2+\alpha)}{\Gamma(n+2+\alpha)}=\frac{n!}{(n+1+\alpha)(n+\alpha)\cdots(2+\alpha)},$$
while $\Gamma(s)$ stands for the usual Gamma function.
In particular, when $\alpha=0$, the (standard weighted) Bergman space $A^2_0(\DD)$
is exactly the (classical) Bergman space $A^2(\DD)$.

From now on, let $N$ be a given positive integer. It is well known that the
shift operator $M_{z^N}$ of multiplicity $N$ on $A^2_\alpha(\DD)(-1<\alpha <+\infty)$ is defined by
$$M_{z^N}f(z)=z^Nf(z),\ \ f\in A^2_\alpha(\DD).$$
When $N=1$, $M_{z^1}$ is often written as $M_z$.

Let $(X,\|\cdot\|_X)$ be a normed linear space and let $M$ is a closed linear subspace
of $X$. Let $B$ be a bounded linear operator on $X$.
If $BM\subset M$, then $M$ is called an invariant
subspace for the operator $B$ on $X$.
If $BM\subset M$ and $B^*M \subset M$, then $M$ is called a reducing subspace for the operator $B$ on X.
If $Y\subset X$, then $[Y]_{B,X}$ denotes the smallest invariant subspace for the operator
$B$ on $X$ which contains the set $Y$. That is,
\begin{eqnarray}\label{e14}
[Y]_{B,X}=\overline{{\mbox{span}}\{B^n y~:~y\in Y, n=0,1,2,\cdot\cdot\cdot\}}^{~\|\cdot\|_X},
\end{eqnarray}
where span$X_0$ denotes the linear span of the subset $X_0$ of the linear space $X$, and
$\overline{X_0}^{~\|\cdot\|_X}$ denotes the norm closure of the subset $X_0$ of the normed linear space $(X,\|\cdot\|_X)$.
In particular, if $(X, \langle \cdot ,\cdot\rangle)$ be a Hilbert space and  $\langle Bx,x\rangle\geq 0$
for every $x\in X$, then $B$ is called a positive operator on $X$.
It is well known that every positive operator $B$ on a complex Hilbert space $X$ is self-adjoint,
and that if a positive operator $B$ is bijective (invertible) from a complex Hilbert space
$X$ onto self, then the inverse operator
$B^{-1}$ is also positive, so it is self-adjoint (cf. \cite{pda04} p.277 and \cite{wr91} p.330-331).

\section{The main result}

\begin{theorem}\label{t01}{
Let $N$ be a given positive integer.
Let $S=M_{z^N}$ be the shift operator of multiplicity $N$ on the Bergman space $A^2_\alpha(\DD)$ $(-1<\alpha < +\infty)$.
If $H$ is a reducing subspace for the operator $S$ on $A^2_\alpha(\DD)$, then
$$H=[H\ominus SH]_{S,A^2_\alpha\left(\DD\right)}.$$
}
\end{theorem}

\textit{\textbf{Proof.}}
Since
$A^2_\alpha(\DD) (-1<\alpha < +\infty)$ is a Hilbert space and since $H$ is a closed linear subspace
of $A^2_\alpha(\DD)$, it follows that $H$ can be regarded as a Hilbert space with the inner product
of $A^2_\alpha(\DD)$. Let $T$ be the restriction of the shift operator $S$ to its invariant subspace $H$, that is,
$$Tf(z)=(S|_{_H} f)(z)=z^N f(z),~~f\in H.$$
Then $TH=SH\subset H$, so that $T$ is a bounded linear operator from the Hilbert space $H$ into self.
We now show that
\begin{eqnarray}\label{e12}
[H\ominus SH]_{S,A^2_\alpha\left(\DD\right)}=[H\ominus TH]_{T,H}.
\end{eqnarray}
In fact, since $SH=TH$, it follows $H\ominus SH=H\ominus TH$. Let
$$E=H\ominus TH.$$
Then $ H\ominus SH=E$. Thus by (\ref{e14}), to show (\ref{e12}) it suffices to show that

\begin{eqnarray}\label{e13}
\nonumber&&\overline{{\mbox{span}}\{S^n f: f\in E, n=0,1,2,\cdots\}}^{~\|\cdot\|_{A^2_\alpha\left(\DD\right)}}\\
\nonumber \\
=&&
\overline{{\mbox{span}}\{T^n f: f\in E, n=0,1,2,\cdots\}}^{~\|\cdot\|_{_H}}.
\end{eqnarray}
Indeed, it is clear that when $f\in E$ we have $f\in H$, and
$T^n f=S^nf \in H$ for each $n=0,1,2,\cdots,$
and that for each $g\in H$ we have
$$\|g\|_{_H}=\|g\|_{A^2_\alpha\left(\DD\right)}.$$
Also, since $H$ is a closed linear subspace of the Bergman
space $A^2_\alpha\left(\DD\right)$, it follows that (\ref{e13}) holds. Consequently, (\ref{e12}) holds.
To show Theorem \ref{t01}, it therefore suffices to prove
\begin{eqnarray}\label{e01}
 H&=& [H\ominus TH]_{T,H}
\end{eqnarray}
(the idea using the operator $T$ is from \cite{asc96} p.286, \cite{pda04} P.276-277, \cite{sms01} p.147-150 and so on).

Moreover, for each $f(z)=\sum^\infty_{n=0}a_nz^n\in H$, we have
\begin{eqnarray}\label{e211}
\nonumber && \left\|Tf\right\|^2_{H}=\left\|S|_Hf\right\|^2=\left\|(M_{z^N}|_H)f\right\|^2_H=\left\|\sum^\infty_{n=0}a_nz^{N+n}\right\|^2_H\\
\nonumber &&=\sum^\infty_{n=0}\omega_{N+n}|a_n|^2= \sum^\infty_{n=0}\frac{N+n}{N+n+1+\alpha}\omega_{N+n-1}|a_n|^2=\cdots\\
&&=\sum^\infty_{n=0}\frac{N+n}{N+n+1+\alpha}\cdot\frac{N+n-1}{N+n+\alpha}\cdot\ldots\cdot\frac{n+1}{n+2+\alpha}\omega_n|a_n|^2
\end{eqnarray}
write
$$C_{N,\alpha,n}=\frac{N+n}{N+n+1+\alpha}\cdot\frac{N+n-1}{N+n+\alpha}\cdot\ldots\cdot\frac{n+1}{n+2+\alpha},~~n=0,1,2,\cdots.$$
It is easy to see that $\frac{1}{3+\alpha}<\frac{m}{m+1+\alpha}<1$ for every $m=1,2,\cdots,$
so that
\begin{eqnarray}\label{e231}
0<\frac{1}{(3+\alpha)^N}<C_{N,\alpha,n}<1.
\end{eqnarray}

Thus by (\ref{e211}), we have

\begin{eqnarray}\label{e02}
\|Tf\|^2_{_H}=
\sum^\infty_{n=0}C_{N,\alpha,n}\omega_n |a_n|^2
\geq \sum^\infty_{n=0}\frac{1}{(3+\alpha)^N}\omega_n |a_n|^2
=\frac{1}{(3+\alpha)^N}\|f\|^2_{_H}.
\end{eqnarray}
for each $f(z)=\sum^\infty_{n=0}a_nz^n\in H$. This implies that the operator $T$ is bounded below
from the Hilbert space $H$ into self. Thus by the bounded inverse theorem (\cite{ckh03} p.75-76), the operator $T$ is injective and has a closed range $TH$. Therefore $TH$
is a closed linear subspace of the Hilbert space $H$.

The proof of Theorem \ref{t01} is divided into three parts.

\textbf{Part 1.} Produce the auxiliary operator $A$.

To this end, we first show that the operator $T^*T$ is a bijection from the Hilbert space $H$ onto self.
We now show that the the operator $T^*T$ is an injection (one-to-one) from $H$ into self. In fact, suppose
$T^*Tf=0$ for some $f\in H$. Then
$$\|Tf\|^2_{_H}=\langle Tf,Tf\rangle=\langle f,T^*Tf\rangle=0.$$
So $Tf=0$. Also, since $T$ is injective, it follows that $f=0$. Consequently
$T^*T$ is injective. Next we show that the operator $T^*T$ is a surjection
from $H$ onto self. In fact, suppose $g\in (T^*TH)^\bot$, that is,
$\langle g,T^*Tf \rangle=0$ for all $f\in H$.
In particular, we have
$$\|Tg\|^2_{_H}=\langle Tg,Tg\rangle=\langle g,T^*Tg \rangle=0.$$
So $Tg=0$.
Again, since $T$ is injective, it follows that $g=0$. Consequently
$(T^*TH)^\perp=\{0\}$.
For convenience, write $M=T^*TH$. Then $M$ is a linear subspace of the Hilbert space $H$.
Therefore, it follows that $M^\perp=\{0\}$ if and only if
$\overline{M}^{~\|\cdot\|_{_H}}=H$. Thus by $(T^*TH)^\bot=\{0\}$,
we have
$\overline{T^*TH}^{~\|\cdot\|_{_H}}=H$. To show that $T^*TH=H$,
it therefore suffices to prove that the operator $T^*T$ has a closed range. Indeed, since
$$\|Tf\|^2_{_H}=\langle Tf, Tf\rangle=\langle f,T^*Tf \rangle \leq \|f\|_{_H}\|T^*Tf\|_{_H}$$
for all $f\in H$, it follows from (\ref{e02}) that
$$\|T^*Tf\|_{_H}\|f\|_{_H}\geq \|Tf\|^2_{_H}\geq \frac{1}{(3+\alpha)^N}\|f\|^2_{_H},$$
so $\|T^*Tf\|_{_H}\geq \frac{1}{(3+\alpha)^N}\|f\|_{_H}$ for all $f \in H$. Therefore the operator $T^*T$ is bounded below
form the Hilbert space $H$ onto self. Hence it has a closed range.
Consequently, $T^*T$ is a bijective bounded linear operator from
the Hilbert space $H$ onto self, that is, it is invertible.
Moreover, it is easy to see that $T^*T$ is a positive operator on $H$, hence its inverse operator
$(T^*T)^{-1}$ is a positive operator on $H$, so $(TT^*)^{-1}$ is a self-adjoint bounded linear operator from the Hilbert space
$H$ onto self.

Now, we can define the auxiliary operator $A$ on the Hilbert space $H$ by
$$Ag=T(T^*T)^{-1}g,~~g\in H.$$
It is clear that $A$ is bounded linear operator on the Hilbert space $H$.
Moreover, since the operator $(T^*T)^{-1}$ is self-adjoint, it follows that
\begin{equation}\label{e03}
 A^*=(T(T^*T)^{-1})^*=((T^*T)^{-1})^*T^*=(T^*T)^{-1}T^*.
\end{equation}

Next, we give an explicit expression for the operator $A$. By the definition,
for any $f(z)=\sum^\infty_{n=0}a_nz^n \in A^2_\alpha(\DD)$ and any
$g(z)=\sum^\infty_{n=0}b_nz^n \in A^2_\alpha(\DD) $, we have
$$M_{z^N}f(z)=\sum^\infty_{n=0}a_nz^{N+n}
=a_{_0}z^N+a_{_1}z^{N+1}+a_{_2}z^{N+2}+\cdots+a_{_n}z^{N+n}+\cdots
,$$
so
\begin{eqnarray*}
  \langle f,M_{z^N}^*g\rangle&=& \langle M_{z^N}f,g\rangle= \sum^\infty_{n=0}\omega_{N+n} a_{n}\overline{b_{N+n}}\\
  &=&\sum^{\infty}_{n=0} \frac{N+n}{N+n+1+\alpha}\cdot\frac{N+n-1}{N+n+\alpha}\cdot\ldots\cdot\frac{n+1}{n+2+\alpha}\omega_na_n\overline{b_{N+n}}.\\
  &=&\sum^{\infty}_{n=0} C_{N,\alpha,n}\omega_na_n\overline{b_{N+n}}
\end{eqnarray*}
This implies that
$$M_{z^N}^*g(z)=\sum^\infty_{n=0}C_{N,\alpha,n}b_{N+n}z^n, ~~g(z)=\sum^\infty_{n=0}b_nz^n\in H.$$
That is
\begin{eqnarray*}
&&M_{z^N}^*(b_0+b_1z+b_2z^2+\cdots+b_nz^n+\cdots) \\
&=&C_{N,\alpha,0}b_{_N}+C_{N,\alpha,1}b_{_{N+1}}z+C_{N,\alpha,2}b_{_{N+2}}z^2+\cdots+
 C_{N,\alpha,n}b_{_{N+n}}z^n+\cdots
\end{eqnarray*}
for any $b_0+b_1z+b_2z^2+\cdots+b_nz^n+\cdots \in A^2_\alpha(\DD)$

Also, since $H$ is a reducing subspace for the operator $S=M_{z^N}$,
it follows from the result of Problem 4.7 in \cite{ckh03} p.36, that
for all $f(z)=\sum_{n=0}^\infty a_nz^n\in H$ and all $g(z)=\sum_{n=0}^\infty b_nz^n\in H$,
we have

\begin{eqnarray*}
&&Tf(z)=(M_{z^N}|_H)f(z)=(S|_H)f(z)=M_{z^N}f(z)\\
&=&a_{_0}z^N+a_{_1}z^{N+1}+a_{_2}z^{N+2}+\cdots+a_{_n}z^{N+n}+\cdots\in H
\end{eqnarray*}
and
\begin{eqnarray*}
&&T^*Tf(z)=(M_{z^N}|_H)^*Tf(z)=(M_{z^N}^*|_H)Tf(z)=M_{z^N}^*Tf(z)\\
&=&M_{z^N}^*(a_0z^N+a_1z^{N+1}+a_2z^{N+2}+\cdots+a_nz^{N+n}+\cdots)\\
&=&C_{N,\alpha,0}a_{_0}+C_{N,\alpha,1}a_{_1}z+C_{N,\alpha,2}a_{_2}z^2+\cdots+C_{N,\alpha,n}a_{_n}z^n+\cdots.
\end{eqnarray*}

Also, since $(T^*T)^{-1}T^*T=I_{_H}$ where $I_{_H}$ denote the identity operator on $H$, It follows that
\begin{eqnarray*}
 &&(T^*T)^{-1}(b_0+b_1z+b_2z^2+\cdots+b_nz^n+\cdots)\\
 &=&\frac{1}{C_{N,\alpha,0}}b_0+\frac{1}{C_{N,\alpha,1}}b_1~z+\frac{1}{C_{N,\alpha,2}}b_2z^2+
 \cdots +\frac{1}{C_{N,\alpha,n}}b_nz^n +\cdots,
\end{eqnarray*}
so
\begin{eqnarray*}
&&Ag(z)=T(T^*T)^{-1}(b_0+b_1z+b_2z^2+\cdots+b_nz^n+\cdots) \\
 &=&\frac{1}{C_{N,\alpha,0}}b_0z^N+\frac{1}{C_{N,\alpha,1}}b_1z^{N+1}+\frac{1}{C_{N,\alpha,2}}b_2z^{N+2}+
 \cdots +\frac{1}{C_{N,\alpha,n}}b_nz^{N+n} +\cdots,
\end{eqnarray*}

\begin{eqnarray*}
&&A^2g(z)=(A(Ag))(z)\\
&=&C_{N,\alpha,0}^{(2)}b_{_0}z^{2N}+C_{N,\alpha,1}^{(2)}b_{_1}z^{2N+1}+C_{N,\alpha,2}^{(2)}b_{_2}z^{2N+2}+
 \cdots +C_{N,\alpha,n}^{(2)}b_{_n}~z^{2N+n} +\cdots,
\end{eqnarray*}

$\cdots\cdots,$

\begin{eqnarray}\label{e04}
A^mg(z)=C_{N,\alpha,0}^{(m)}b_0z^{mN}+C_{N,\alpha,1}^{(m)}b_1z^{mN+1}+
 \cdots +C_{N,\alpha,n}^{(m)}b_nz^{mN+n} +\cdots,
\end{eqnarray}

$\cdots\cdots,$\\
where $g(z)=\sum^\infty_{n=0}b_nz^n\in H$, $m=2,3,\cdots$, while $C_{N,\alpha,0}^{(m)}, C_{N,\alpha,1}^{(m)}, C_{N,\alpha,2}^{(m)},
 \cdots, C_{N,\alpha,n}^{(m)}$ are some constants greater than 1 which are dependent on
$N$, $\alpha$ and $m$.

\textbf{Part 2.} Properties of the auxiliary operator $A$.

Property 1. $A^nH$ is a closed linear subspace of the Hilbert space $H$ for each $n=1,2,\cdots.$

In fact, by the above formula, we have
$$Ag(z)=\sum^\infty_{n=0}\frac{1}{C_{N,\alpha,n}}b_nz^{N+n} ,~~~g(z)=\sum^\infty_{n=0}b_nz^n\in H.$$
Thus by (\ref{e211}) and (\ref{e231}), we have
$$\|Ag\|^2_{_H}=\sum^\infty_{n=0}\omega_{N+n}\left(\frac{1}{C_{N,\alpha,n}}|b_n|\right)^2=\sum^\infty_{n=0}\frac{1}{C_{N,\alpha,n}}\omega_n|b_n|^2
\geq\sum^\infty_{n=0}\omega_n|b_n|^2=\|g\|^2_{_H}$$
for every $g\in H$. Hence
$$\|A^2g\|_{_H}=\|A(Ag)\|_{_H}\geq \|Ag\|_{_H}\geq\|g\|_{_H}, ~~g\in H.$$
In general, we have
$$\|A^ng\|_{_H}\geq \|g\|_{_H}, ~~g\in H, n=1,2,\cdots.$$
That is, the operator $A^n$ is bounded below from the Hilbert space $H$ into self
for each $n=1,2,\cdots$. Therefore $A^nH$ is a closed linear subspace of the Hilbert
space $H$ for each $n=1,2,\cdots$.

Property 2. Let
$E=H\ominus TH$. Then
\begin{eqnarray}\label{e06}
{\mbox{ker}}(A^*)^n\subset E+TE+\cdots +T^{n-1}E, ~~~n=1,2,\cdots.
\end{eqnarray}

To prove (\ref{e06}), we first show that $E=\mbox{ker}T^*$, which can be obtained by
$g\in E=H\ominus TH~\Leftrightarrow~\langle Tf,g\rangle=0$ for all $f\in H~\Leftrightarrow~\langle f,T^*g\rangle=0$
for all $f\in H$ $\Leftrightarrow$ $T^*g=0$ $\Leftrightarrow$ $g\in \mbox{ker}T^*$.

Next, by (\ref{e03}) we have
\begin{eqnarray}\label{e07}
A^*T=(T^*T)^{-1}T^*T=I_{_H},
\end{eqnarray}
where $I_{_H}$ denote the identity operator on the Hilbert space $H$.
For each $g\in TH$, there is an $f\in H$ such that $g=Tf$. Thus by (\ref{e07})
we have
\begin{eqnarray}\label{e08}
TA^*g=TA^*Tf=Tf=g,~~g\in TH.
\end{eqnarray}
On the other hand, since $E=\mbox{ker}T^*$, it follows that $T^*E=\{0\}$. Thus by
(\ref{e03}) we obtain
\begin{eqnarray}\label{e09}
TA^*g=T(T^*T)^{-1}T^*g=0,~~~g\in E.
\end{eqnarray}
Since $TH$ is a closed linear subspace of the Hilbert space $H$ and since
$H=TH\oplus E$, it follows from (\ref{e08}) and (\ref{e09})
that $P_{_{TH}}=TA^*$ and
$$P_E=I_{_H}-P_{_{TH}}=I_{_H}-TA^*,$$
where $P_{_{TH}}$ and $P_E$ denote the orthogonal projection operators from the Hilbert space $H$ onto its closed linear subspaces $TH$ and $E$ respectively.
Consequently, we have
$$T^kP_E(A^*)^k=T^k(A^*)^k-T^{k+1}(A^*)^{k+1}, ~~~k=0,1,2\cdots.$$
Therefore
\begin{eqnarray*}
&&\sum^{n-1}_{k=0}T^kP_E(A^*)^k\\
&=&[I_{_H}-TA^*]+[TA^*-T^2(A^*)^2]+\cdots +[T^{n-1}(A^*)^{n-1}-T^n(A^*)^n]\\
&=&I_{_H}-T^n(A^*)^n.
\end{eqnarray*}
Now, let $g\in \mbox{ker}(A^*)^n$. Then $(A^*)^ng=0$. Noting that $P_E(A^*)^kg \in E~(k=0,1,2,\cdots),$ we have
$$g=[I_{_H}-T^n(A^*)^n]g=\sum^{n-1}_{k=0}T^kP_E(A^*)^kg \in \sum^{n-1}_{k=0}T^kE=E+TE+\cdots+T^{n-1}E.$$
That is, $\mbox{ker}(A^*)^n\subset E+TE+\cdots +T^{n-1}E$ for each $n=1,2,\cdots$.

Property 3.
\begin{eqnarray*}\label{e05}
\bigcap_{n=1}^\infty A^nH=\{0\}.
\end{eqnarray*}

In fact, if not, then there is at least a nonzero function $f_0$ in $H$ such that
$$f_0\in A^n H, ~~n=1,2,\cdots.$$
This implies that $f_0(z)=A^ng_n(z)$ for some $g_n\in H$ and every $n=1,2,\cdots$.
Thus by (\ref{e04}), we have
$$f_0(z)=z^{nN} h_n(z),~~n=1,2,\cdots,$$
where $h_n(z)$ is an analytic function in the open unit disk $\DD$ on the complex
plane. Thus the analytic function $f_0(z)$ has a zero of $nN$ order at the origin
for each $n=1,2,\cdots$, which contradicts the fact that any nonzero analytic
function can have at most a zero of finite order at the origin (\cite{wr87} p.208-209, Theorem 10.18).

\textbf{Part 3}. Proof of the expression (\ref{e01}).

By Property 1, $A^nH$ is a closed linear subspace of the Hilbert space $H$ for each $n=1,2,\cdots$. Thus by the
projection theorem, for each $f\in H$ and each $n=1,2,\cdots$, we can uniquely write
\begin{eqnarray}\label{e10}
f=g_n+h_n,
\end{eqnarray}
where $g_n\in A^nH$ and $h_n\in (A^nH)^\bot$. Thus by Property 2,
\begin{eqnarray}\label{e11}
h_n\in (A^nH)^\bot=\mbox{ker}(A^*)^n\subset E+TE+\cdots +T^{n-1}E\subset [E]_{T,H},
\end{eqnarray}
where $E=H\ominus TH$, while $[E]_{T,H}$ denotes the smallest invariant subspace for the operator $T$ on $H$ which
contains the set $E$. Moreover, by (\ref{e10}) we have
$$\|f\|^2_{_H}=\|g_n\|^2_{_H}+\|h_n\|^2_{_H}, ~~~n=1,2,\cdots.$$
So $\|g_n\|_{_H}\leq \|f\|_{_H}$ for each $n=1,2,\cdots$, that is, $\{g_n\}$ is a bounded sequence
in the Hilbert space $H$. Hence there is a subsequence $\{g_{n_k}\}$ of $\{g_n\}$ that
converges weakly to some element $g\in H$, that is
\begin{eqnarray}\label{e15}
\langle g_{n_k},w\rangle\rightarrow \langle g,w\rangle
\end{eqnarray}
for every $w\in H$ as $k\rightarrow \infty$.

We now show that $g=0$. To this end, suppose on the contrary that $g\neq 0$. But
by Property 3,
$\bigcap_{n=1}^\infty A^nH=\{0\}$, so
$$g\notin \bigcap_{n=1}^\infty A^nH.$$
This implies that there is at least a positive integer $m$ such that
$g\notin A^mH$. Again by the projection theorem, $g$ can be uniquely written
as
$$g=u+v,$$
where $u\in A^mH$ and $v\in (A^mH)^\perp$. Since $g\notin A^mH$, it follows that $v\neq 0$
(otherwise, $g=u\in A^mH$). Therefore, we have
\begin{eqnarray}\label{e16}
\langle g,v\rangle=\langle u+v,v\rangle=\langle v, v\rangle=\|v\|^2_{_H}\neq 0.
\end{eqnarray}
On the  other hand, $g_{n_k}\in A^{n_k}H\subset A^mH$ for all $n_k\geq m$, so $\langle g_{n_k},v\rangle=0$
for all $n_k\geq m$. Thus by (\ref{e15}) and (\ref{e16}), we have
$$0=\lim_{k\rightarrow \infty}\langle g_{n_k},v\rangle=\langle g,v\rangle\neq 0.$$
This is a contradiction. Therefore we conclude that the weak limit $g=0$.

Moreover, it follow from the basic theory of Banach spaces (see, for example, \cite{rem98} p.217, Corollary 2.5.19)
that if $\{x_n\}$ is a sequence in a normed space $X$ that converges weakly to some $x\in X$, then there is a sequence
$\{y_n\}$ of convex combinations of members of $\{x_n:n=1,2,\cdots\}$ that converges to $x$ in the norm of $X$.
Now, since the sequence $\{g_{n_k}\}$ in the Hilbert space $H$ converges weakly to $g\in H$,
there is a sequence $\{u_k\}$ of convex combinations of members
of $\{g_{n_k}:k=1,2,\cdots\}$ that converges to $g$ in the norm of $H$. Thus we can assume that
$$u_k=c_{k_1}g_{n_{k_1}}+c_{k_2}g_{n_{k_2}}+ \cdots +c_{k_l}g_{n_{k_l}}\rightarrow g=0$$
as $k\rightarrow \infty$ in the norm of $H$, where $c_{k_1}$, $c_{k_2},\cdots,$ $c_{k_l}$
are nonnegative real numbers with $c_{k_1}+c_{k_2}+\cdots+c_{k_l}=1$, while
$g_{n_{k_1}}$, $g_{n_{k_2}},\cdots,~g_{n_{k_l}}$ are some members of $\{g_{n_{k}}:k=1,2,\cdots\}$.

On the other hand, by (\ref{e10}) and (\ref{e11}), we have
$$f=g_{n_{k_j}}+h_{n_{k_j}},~~~j=1,2,\cdots,l,$$
where $g_{n_{k_j}}\in A^{n_{k_{j}}} H$ and $h_{n_{k_j}}\in (A^{n_{k_j}}H)^\bot\subset [E]_{T,H}$.
Write
$$v_k=c_{k_1}h_{n_{k_1}}+c_{k_2}h_{n_{k_2}}+ \cdots +c_{k_l}h_{n_{k_l}}.$$
Then $v_k \in [E]_{T,H}$ and
\begin{eqnarray*}
f&=&(c_{k_1}+c_{k_2}+\cdots+c_{k_l})f  \\
&=&c_{k_1}(g_{n_{k_1}}+h_{n_{k_1}})+ c_{k_2}(g_{n_{k_2}}+h_{n_{k_2}})+\cdots+c_{k_l}(g_{n_{k_l}}+h_{n_{k_l}})=u_k+v_k.
\end{eqnarray*}
Also, since $u_k\rightarrow 0$ as $k\rightarrow \infty$ in the norm of $H$, it follows that
$$f=\lim_{k\rightarrow \infty}(f-u_k)=\lim_{k\rightarrow \infty}v_k\in [E]_{T,H},$$
where the limit is taken in the norm of $H$ (i.e. in the norm of $A^2_{\alpha}~(\DD)(-1<\alpha<+\infty))$.
That is,
$H\subset [E]_{T,H}$.
Moreover, since $E=H\ominus TH$, it follows from (\ref{e14}) that
$$[E]_{T,H}=\overline{{\mbox{span}}\{T^n f: f\in E, n=0,1,2,\cdots\}}^{~\|\cdot\|_{_H}}\subset H.$$
Thus by (\ref{e12}), we have
$$H=[E]_{T,H}=[H\ominus TH]_{T,H}=[H\ominus SH]_{S,A^2_\alpha\left(\DD\right)}.$$
The proof of the theorem is finished.

\end{document}